\documentclass[12pt] {article}
\usepackage{amsmath,amsthm}
\usepackage{comment}
\title{Range characterization of the cosine transform on higher Grassmannians. }
\date{}
\author{ Semyon Alesker , Joseph Bernstein
\\  { \normalsize Department of Mathematics, Tel Aviv University, Ramat Aviv}
 \\  { \normalsize 69978 Tel Aviv,
Israel }}
\newcommand{\RR}{\mbox{\rm $~\vrule height6.5pt width0.5pt
depth0.3pt\!\!$R}}

\newcommand{\CC}{\mbox{\rm $~\vrule height6.5pt width0.5pt
depth0.3pt\!\!$C}}

\def\lam{\lambda}

\def\str{\longrightarrow}

\def\qed { Q.E.D. }

\def\gr{ Gr}
\def\val{Val^{ev}_{n,i}}
\swapnumbers
\newtheorem{theorem}{Theorem}[section]

\newtheorem{lemma}[theorem]{Lemma}
\newtheorem{proposition}[theorem]{Proposition}

\theoremstyle{definition}

\newtheorem{definition}[theorem]{Definition}

\begin{document}
\maketitle
\begin{abstract}
We characterize the range of the cosine transform on real
Grassmannians in terms of the decomposition under the action of
the special orthogonal group $SO(n)$. We also give a geometric
interpretation of this image in terms of valuations. In addition,
we discuss the non-Archimedean analogues.
\end{abstract}
\setcounter{section}{-1}
\section{Introduction.}
The main results of this paper are Theorems 1.2, 1.3, and 2.1
below. In Section 1 we study the image of the cosine transform.
The investigations of this problem were started by G. Matheron in
1974 \cite{matheron1}. It was investigated further by Goodey and
Howard \cite{goodey-howard1}, \cite{goodey-howard2}, Goodey,
Howard, and Reeder \cite{goodey-howard-reeder}. For more details
on the previous results  see Section 1 of this paper. Theorem 1.2
describes the image of the cosine transform (defined in Section 1)
acting from the space of functions on the Grassmannian of real
$j$- dimensional subspaces in $\RR^n$ to the space of functions on
the Grassmannian of real $i$-dimensional subspaces. The
description is given in terms of $K$-types, namely in terms of the
decomposition into irreducible subspaces under the action of the
special orthogonal group $SO(n)$. In Theorem 1.3 we consider
separately the case $i=j$; this case is especially important since
the study of other cases is reduced to this one using the results
on the Radon transform due to Gelfand, Graev, and Ro\c su
\cite{gelfand-graev-rosu} (see also Grinberg \cite{grinberg}). It
turns out that the cosine transform in this case can be
interpreted as an intertwining operator of $GL(n,\RR)$-modules.
This key observation allows us to use the representation theory of
reductive groups.
 We prove that this intertwining operator has an irreducible image.
 This step uses connection of the cosine transform with the theory
 of valuations on convex sets (discussed in Section 1) and earlier
 results of one of the authors about representation theoretical properties of valuations \cite{alesker1},
 \cite{alesker}. Once we have this irreducibility result the
 description of the $K$-type structure is a direct consequence
 of the results of Howe and Lee \cite{howe-lee}.

In Section 2 we discuss an analogue of the cosine transform over
non-Archimedean local fields.

The cosine transform
 appears naturally in convex,
integral, and stochastic geometry (see \cite{goodey-howard1},
\cite{goodey-howard2}, \cite{goodey-howard-reeder},
\cite{goodey-zhang}, \cite{groemer},\cite{matheron1},
\cite{matheron2}, \cite{schneider70}). For the basic notions of
the representation theory  we refer to \cite{wallach}.

{\bf Remarks.} (a) The results of this paper were applied further
to the theory of valuations in \cite{alesker-univa}.

(b) Recently we have found a more general result about the cosine
type transform of more general form (where the kernel $|cos(E,F)|$
is replaced by $|cos (E,F)|^{\alpha},\, \alpha \in \CC$; the
notation is explained in Section 1). This will be discussed
elsewhere.

{\bf Acknowledgements.} The first named author is grateful to A.
Braverman, D. Kazhdan, and D. Soudry for very useful discussions.
We thank E. Grinberg for useful remarks.

\section{Range of the cosine transform.}
First let us recall some notation. We will denote by $\gr_{k,n}$
the Grassmannian of real $k$- dimensional subspaces in the real
$n$- dimensional space $\RR^n$. Let us fix a Euclidean structure
on $\RR^n$. Let $E\in \gr_{i,n}, \, F\in \gr_{j,n}$. Assume that
$i\leq j$. Let us call by {\itshape cosine of the angle} between
$E$ and $F$ the following number:
$$|cos(E,F)|:=\frac{vol_i(Pr_F(A))}{vol_i(A)},$$
where $A$ is any subset of $E$ of non-zero volume, $Pr_F$ denotes
the orthogonal projection onto $F$, and $vol_i$ is the $i$-
dimensional measure induced by the Euclidean structure. (Note that
this definition does not depend on the choice of a subset
$A\subset E$). In the case $i\geq j$ we define the cosine of the
angle between them as cosine of the angle between their orthogonal
complements:
$$|cos(E,F)| :=|cos(E^\perp ,F^\perp )|.$$
(It is easy to see that if $i=j$ both definitions are equivalent.)

Let us call by {\itshape  sine of the angle} between $E$ and $F$
the cosine between $E$ and the orthogonal complement of $F$:
$$|sin(E,F)|:= |cos(E,F^\perp )|.$$
 The following properties are well known (and rather trivial):
$$|cos(E,F)|=|cos(F,E)|=|cos(E^\perp ,F^\perp)|,$$
$$|sin(E,F)|=|sin(F,E)|=|sin(E^\perp ,F^\perp)|,$$
$$0\leq |cos(E,F)|, \, |sin(E,F)| \leq 1.$$

For any $1\leq i, \, j \leq n-1$ one defines the cosine transform
$$T_{j,i}:C(\gr_{i,n}) \str C(\gr_{j,n})$$ as follows:
$$(T_{j,i}f)(E):= \int_{\gr_{i,n}} |cos(E,F)| f(F) dF,$$
where the integration is with respect to the Haar measure on the
Grassmannian.

{\bf Remark.} One should notice that very often in the literature
the cosine transform $T_{j,i}$ (resp. the Radon transform
$R_{j,i}$) is denoted by $T_{i,j}$ (resp. $R_{i,j}$), i.e. with
permutation of indexes. We prefer our notation since it is more
convenient to write the composition formulas like
$R_{j,i}=R_{j,k}\circ R_{k,i}$.

 Clearly the cosine transform  commutes with the
action of the orthogonal group  $O(n)$,  and hence its image is
$O(n)$- invariant subspace of functions. Our first main result
(Theorem 1.2) is the description of the image of the cosine
transform in terms of the decomposition of it with respect to the
action of $SO(n)$. Since the representation of $SO(n)$ in
functions on the Grassmannian is multiplicity free, it is
sufficient to list those irreducible representations of $SO(n)$
entering into the image of the cosine transform. Moreover it is
shown (Theorem 1.3) that for $i=j$ this image coincides with the
image of even translation invariant $i$- homogeneous continuous
valuations on convex sets under certain natural map (actually this
fact is proved first and then used in the proof of the more
explicit result).

Now let us recall standard facts on the representations of the
special orthogonal group $SO(n)$ (see e.g. \cite{zhelobenko}).
\begin{lemma}
 The isomorphism classes of irreducible representations of
$SO(n)$,
$n>2$ are parameterized by their highest weights, namely
sequences of integers $(m_1, m_2, \dots ,m_{[n/2]})$ which
satisfy:

(i) if $n$ is odd then
$m_1 \geq m_2  \geq \dots \geq  m_{[n/2]} \geq 0$:

(ii) if $n>2$ is even then
$m_1 \geq m_2 \geq \dots \geq  m_{n/2 -1} \geq |m_{n/2}|$.
\end{lemma}
Recall also that for $n=2$ the representations of $SO(2)$ are
parameterized by a single integer $m_1$. We will use the following
notation. Let us denote by $\Lambda^+$ the set of all highest
weights of $SO(n)$, and by $\Lambda^+_k$ the set of all highest
weights $\lambda =(m_1, m_2, \dots , m_{[n/2]})$ with $m_i=0$ for
$i>k$ and all $m_i$ are {\itshape even}.

Let us recall the decomposition of the space of functions on the
Grassmannian $\gr_{k,n}$ under the action of $SO(n)$ referring for
the proofs to \cite{sugiura}, \cite{takeuchi}. Since $\gr_{k,n}$
is a symmetric space, each irreducible representation enters with
multiplicity at most one. The representations which do appear have
highest weights precisely from $\Lambda^+_k \cap \Lambda^+_{n-k}$.
\newline
Now let us state our main result.
\begin{theorem}
Let $1\leq i, j \leq n-1$. Then the image of the cosine transform
$T_{j,i} :C(\gr_{i,n} )\str C(\gr_{j,n})$ consists of irreducible
representations of $SO(n)$ with highest weights $\lambda =(m_1,
\dots , m_{[n/2]})$ precisely with the following additional
conditions:

(i)$\lambda \in \Lambda^+_{i} \cap \Lambda^+_{n-i} \cap
\Lambda^+_{j} \cap \Lambda^+_{n-j}$;

(ii) $|m_2| \leq 2$;

Moreover every $C^\infty$- function on $\gr_{j,n}$ which belongs
to the closure of the sum of $SO(n)$- irreducible subspaces
satisfying the above conditions (i)-(ii), is an image under
$T_{j,i}$ of some $C^\infty$- function on $\gr_{i,n} $.
\end{theorem}
Note that as a corollary of this theorem we immediately get
exactly the same characterization of the image of the sine
transform. Theorem 1.2 was known for a long time for $i=j= 1$ (or
equivalently for $i=j=n-1$), see \cite{schneider70}
 or  \cite{groemer}. The case $n=4, \, i=j=2$ was described completely in
\cite{goodey-howard-reeder}; this paper contains also partial
information on the general case.

The main case in the proof of this theorem is $i=j$. The cosine
transform for $i \ne j$ can be written as a composition of the
Radon transform between different Grassmannians and the cosine
transform for $i=j$. Thus to deduce the general case we use the
characterization of the image of the Radon transform
\cite{gelfand-graev-rosu}, see also \cite{grinberg}.  In order to
treat the case $i=j$ we first interpret the cosine transform as an
intertwining operator between certain representations on the
(non-compact) group $GL(n,\RR)$ induced from maximal parabolic
subgroups. Next we prove that the image of this operator is an
irreducible $GL(n,\RR)$- module . In order to do it we first show
(see Theorem 1.3 ) that the image is contained in the subspace
corresponding to even translation invariant valuations (see the
definitions below). But the last space is an irreducible
$GL(n,\RR)$- module  by the main result in \cite{alesker}. Hence
it coincides with the image of the cosine transform. The
decomposition of the space of even valuations with respect to the
action of $SO(n)$ was described in \cite{alesker} as an easy
corollary of the irreducibility and the computations due to Howe
and Lee \cite{howe-lee}.

Now let us describe the construction of the intertwining operator.
Let us denote by $L$ the line bundle over the Grassmannian $\gr_{i,n}$
whose fiber over a subspace $E\in \gr_{i,n}$ is the space of Lebesgue
measures on $E$ (which is denoted by $|\wedge ^i E^*|$). Clearly $L$ is $GL(n,\RR)$- equivariant
line bundle over $\gr_{i,n}$. Moreover if we fix the Euclidean structure
on $\RR^n$ we can identify $L$ with the trivial bundle in a way compatible with
the action of $SO(n)$. Let $M$ denote the line bundle
over the Grassmannian $\gr_{n-i,n}$ whose fiber over
$F\in  \gr_{n-i,n}$ is the space of Lebesgue measures on
the quotient space $\RR^n /F$ denoted by $|\wedge ^{i} (\RR^n /F)^*|$.
 Let $|\omega|$ denote the line bundle
of densities over $\gr_{n-i,n}$. Let $N:=L\otimes |\omega|$.
Define an intertwining operator $T$ from the space of continuous
sections of $N$ to the space of continuous sections of $L$
$$T:\Gamma (\gr_{n-i,n},N) \str \Gamma (\gr_{i,n}, L) $$
as follows. For $E\in \gr_{i,n}$ and $f\in \Gamma (\gr_{n-i,n},N)$ set
$$(Tf)(E)=\int _{F\in \gr_{n-i,n}} pr^*_{E,F} (f(F)) ,$$
where $pr_{E,F}$ denotes the natural map $E\str \RR^n /F$ and
$pr^*_{E,F}$ is the induced map $|\wedge ^{i} (\RR^n /F)^*| \str |\wedge ^i E^*|$.
Clearly $T$ is a non-trivial operator commuting with the action of $GL(n,\RR)$.
\begin{theorem}
The image of the operator $T$ is an irreducible $GL(n,\RR)$- module.
Moreover if we identify $L$ with the trivial bundle in an
$SO(n)$- equivariant way then the image of $T$ coincides
with the image of the cosine transform $T_{i,i}$.
\end{theorem}
Note that the second statement of the theorem easily follows
from the definitions. In order to prove the irreducibility of the image
we will need one more construction which we are going to
describe.

Let ${\cal K}^n$ denote the family of all convex compact subsets in $\RR^n$.
\begin{definition}
1) A function $\phi :{\cal K}^n \str \CC$ is called a valuation if
for any $K_1, \, K_2 \in {\cal K}^n$ such that their union is also
convex one has
$$\phi(K_1 \cup K_2)= \phi(K_1) +\phi(K_2) -\phi(K_1 \cap K_2).$$

2) A valuation $\phi$ is called continuous if it is continuous with
respect the Hausdorff metric on ${\cal K}^n$.

3) A valuation $\phi$ is called translation invariant if
 $\phi (K+x)=\phi(K)$ for every $x\in \RR^n$ and every $K$.

4) A valuation $\phi$ is called even if $\phi(-K)=\phi(K)$
for every $K\in {\cal K}^n$.

5) A valuation $\phi$ is called homogeneous of degree $k$
(or k-homogeneous) if for every $K\in {\cal K}^n$ and
every scalar $\lam \geq 0$
$\phi(\lam \cdot K)=\lam ^k \phi(K)$.
\end{definition}

We refer to further details on valuations to the surveys
\cite{mcmullen-schneider} and \cite{mcmullensurvey}.
We will need only few facts about them.
\begin{lemma}\cite{hadwiger}
Every translation invariant $n$- homogeneous continuous
valuation on ${\cal K}^n$ is proportional to the Lebesgue measure on $\RR^n$.
\end{lemma}
Let us denote by $Val^{ev}_{n,i}$ the linear space of translation
invariant $i$- homogeneous even continuous valuations. It is a Frechet
space with respect to the topology of uniform
convergence on compact subsets of ${\cal K} ^n$.

There is a natural map
$$\gamma :\val \str \Gamma(\gr_{i,n}, L),$$
where $L$ is the line bundle defined  above. To define it, fix any
$\phi \in \val$. Take any $E\in \gr_{i,n}$. Consider the
restriction of $\phi$ to the class of convex compact subsets of
$E$. It is again translation invariant $i$- homogeneous continuous
valuation on $E$. Hence by Lemma 1.5 it is proportional to the
Lebesgue measure on $E$. Hence $\phi$ defines a continuous section
$\gamma (\phi)$ of $L$. This map was used in \cite{alesker1} and
\cite{alesker} in the proof of McMullen's conjecture. This map was
also independently considered by D. Klain \cite{klain2}. Note that
this map turns out to be injective by a theorem of D. Klain
\cite{klain1}. The main fact on $\val$ we use is the following
result proved in \cite{alesker}.
\begin{lemma}
For every integer $i$, $1\leq i \leq n-1$ the space $\val$ is an irreducible
$GL(n,\RR^n)$- module. Hence its image in $\Gamma (\gr_{i,n}, L)$
is an irreducible submodule.
\end{lemma}

{\bf Proof} of Theorem 1.3. By Lemma 1.6 it remains to show that
the image of $T$ is contained in $\gamma (\val)$. Fix any $f\in
\Gamma (\gr_{n-i,n}, N)$. Let us define a valuation $\phi$ as
follows. For every $K\in {\cal K}^n$ set
$$\phi (K):= \int _{F\in \gr_{n-i,n}} f(pr_{\RR^n /F} (K)),$$
where $pr_{\RR^n /F} :\RR^n \str \RR^n /F$ is the canonical map.
It is easy to see that $\phi \in \val$; and moreover $\gamma(\phi)
=T(f)$. Thus Theorem 1.3 is proved. \qed

{\bf Proof} of Theorem 1.2. First consider the case $i=j$. By
Theorem 1.3 and its proof the image of the cosine transform
coincides with the image under the map $\gamma$ of translation
invariant $i$- homogeneous continuous valuations. The explicit
decomposition of the last space under the action of $SO(n)$ was
given in \cite{alesker} (it was heavily based on \cite{howe-lee}).
Thus Theorem 1.2 is proved for $i=j$.
Now consider the case $i \ne j$. Clearly we may assume that $j<i$.
One has the Radon transform
$$R_{j,i}: C(\gr_{i,n})\str C(\gr_{j,n})$$
defined by $(R_{j,i}f)(H)= \int_{F \supset H}  f(F) \cdot dF$. The
next lemma is well known, but we will present a proof for
convenience of the reader.
\begin{lemma}
Let $1\leq j<i \leq n$. Then
$$T_{j,i}= c \cdot T_{j,j} R_{j,i},$$
here $c$ is a constant depending on $n,i,j$.
\end{lemma}
{\bf Proof.} Fix a subspace $E\in \gr_{j,n}$. Let $A$ be any
convex compact subset of $E$ of positive measure. Let $f$ be a
continuous function on $\gr_{i,n}$. Then by definition
\begin{equation}
 (T_{j,i}f)(E)=\int_{F\in \gr_{i,n}} \frac{vol_j (Pr_F(A))}{vol_j (A)} f(F) dF.
\end{equation}
Let us fix $F\in \gr_{i,n}$ for a moment. Let $B:= Pr_F(A)$.
 By the Cauchy- Kubota formula (see e.g. \cite{schneider-weil})
$$vol_j (B) =c \cdot \int _{H\in Gr_j(F)} vol_j (Pr_H (B) )dH,$$
where $ Gr_j(F)$ denotes the Grassmannian of $j$- dimensional
subspaces of $F$, and $c$ is a non-zero normalizing constant
depending on $i,j$, and $n$ only. Substituting this into (1) and
using the standard change of integration (normalizing all the Haar
measures to be probability measures) one gets:
$$  (T_{j,i}f)(E)=\frac{c}{vol_j(A)} \int_{F\in \gr_{i,n}} dF \cdot f(F) \left(
\int_{H\subset F} dH \cdot vol_j (Pr_F (A)) \right)=$$
$$ \frac{c}{vol_j(A)} \int_{H\in \gr_{j,n}} dH \cdot vol_j (Pr_H (A) )
\left( \int _{F\supset H} dF \cdot f(F) \right).$$ This identity
clearly proves the lemma. \qed

To finish the proof of Theorem 1.2 we need the following fact
proved in
  \cite{gelfand-graev-rosu} (see also \cite{grinberg}).

\begin{proposition}
For $j<i$ the Radon transform $R_{j,i}: C(\gr_{i,n})\str
C(\gr_{j,n})$ is injective iff $i+j \geq n$ and has a dense image
iff $i+j \leq n$.
\end{proposition}

This proposition, description of the decomposition of the space of
functions on the Grassmannians under the action of $SO(n)$, and
the characterization of the image of the cosine transform for
$i=j$ imply the first part of Theorem 1.2.

It remains to prove that if $f$ is a $C^\infty$- function on
$\gr_{i,n}$ belonging to the closure of the sum of $SO(n)$-
irreducible subspaces satisfying the conditions (i)-(ii) in the
statement of Theorem 1.2, then $f$
 is an image under $T_{j,i}$ of some $C^\infty$- function
on $\gr_{j,n} $. We may assume that $j<i$. By Lemma 1.7 $T_{j,i}=c
\cdot T_{j,j} R_{j,i}$. We will need the following fact due to
Casselman and Wallach (see \cite{casselman}).
\begin{proposition} Let $G$ be a real reductive group. Let $K$ be its maximal
compact subgroup. Let $\xi: X \str Y$ be a morphism of two
admissible Banach $G$- modules of finite length which has  a dense
image. Then $\xi$ induces an epimorphism on the spaces of smooth
vectors.
\end{proposition}

In our situation we will need the following more precise form of
Proposition 1.9. In fact it can be proved in a much more general
context, but we do not need it and do not have precise reference.

\begin{lemma}
Let $G= GL(n, \RR)$. Let $K=O(n)$ be the maximal compact subgroup.
Let $X$ and $Y$ be $G$- modules of continuous sections of some
finite dimensional $G$- equivariant vector bundles over the
Grassmannians (or any other partial flag manifolds). Let $\xi: X
\str Y$ be a morphism of these $G$- modules. Then if $f\in Y$ is a
smooth vector belonging to the closure in the continuous topology of the image of $\xi$, 
  then there exists a smooth vector $g\in X$ such
that $\xi(g)=f$ and the $K$- types entering into the decomposition
of $g$ are the same as those of $f$.
\end{lemma}

We postpone the proof of this lemma till the end of the section.
Now let us continue the proof of Theorem 1.2. We have given an
interpretation of $T_{j,j}$ as an intertwining operator of two
$GL(n,\RR)$- modules; they satisfy the conditions of Proposition
1.10 , since they are induced from characters of parabolic
subgroups (see \cite{wallach}). Hence by Proposition 1.10 there
exists a $C^\infty$- smooth function $g$ on the Grassmannian
$Gr_{j,n}$ such that $f= T_{j,j} (g)$ and with the same $K$- types
as $f$.
 Next there exists an interpretation of the Radon
transform as an intertwining operator of some admissible
$GL(n,\RR)$- modules of finite length (it was given in
\cite{gelfand-graev-rosu}). Hence Proposition 1.10 implies the
statement. \qed

{\bf Proof} of Lemma 1.10. Let $f\in Y$ be as in the statement of
Proposition 1.10. By Proposition 1.9 we can choose a smooth vector
$h\in X$ such that $\xi (h)=f$. Then $h$ is just a smooth section
of the corresponding vector bundle over the Grassmannian (or the
partial flags manifold). Let us consider the image $g$ under the
orthogonal projection of $h$ to the closure with respect to the
$L_2$- metric of the span in $X$ of those $K$- irreducible
subspaces which have the same $K$- types as those entering into
the decomposition of $f$. Clearly $\xi(g)=\xi(h)=f$. To finish the
proof it remains to show that $g$ is also a smooth section. Indeed
since the orthogonal projection commutes with the action of $K$,
$K$- smooth vectors go to $K$- smooth, hence $g$ is $K$- smooth
vector. But since $K$ acts transitively on the Grassmannians (and
in fact on all partial flag manifolds), every $K$- smooth section
of every $K$- equivariant bundle must be smooth in the usual
sense. \qed

\section{Non- Archimedean analogue of the cosine transform.}
In this section we study non- Archimedean analogue of the cosine
transform. More precisely we study a non-Archimedean  analogue of
the intertwining operator $T_{i,i}$ (in the notation of the
previous section). We show that it has an irreducible image.

Now let us introduce the necessary notation. Let $F$ be a non-
Archimedean local field. In this section we will denote by
$Gr_{i,n}$ the Grassmannian of $i$- dimensional subspaces in
$F^n$. Denote by $L$ the line bundle over the Grassmannian
$\gr_{i,n}$ whose fiber over a subspace $E\in \gr_{i,n}$ is the
space of Lebesgue measures on $E$ (which is denoted by $|\wedge ^i
E^*|$). Clearly $L$ is $GL(n, F)$- equivariant line bundle over
$\gr_{i,n}$. Let $M$ denote the line bundle over the Grassmannian
$\gr_{n-i,n}$ whose fiber over $H \in  \gr_{n-i,n}$ is the space
of Lebesgue measures on the quotient space $F^n /H$ denoted by
$|\wedge ^{i} (F^n /H)^*|$.
 Let $|\omega|$ denote the line bundle
of densities over $\gr_{n-i,n}$. Let $N:=M \otimes |\omega|$.
Define an intertwining operator $T$ from the space of smooth (i.e., locally constant)
sections of $N$ to the space of continuous sections of $L$
$$T:\Gamma (\gr_{n-i,n},N) \str \Gamma (\gr_{i,n}, L) $$
as follows. For $E\in \gr_{i,n}$ and $f\in \Gamma (\gr_{n-i,n},N)$ set
$$(Tf)(E)=\int _{H\in \gr_{n-i,n}} pr^*_{E,H} (f(H)) ,$$
where $pr_{E,H}$ denotes the natural map $E\str F^n /H$ and
$pr^*_{E,H}$ is the induced map $|\wedge ^{i} (F^n /H)^*| \str
|\wedge ^i E^*|$. Clearly, $T$ is a non-trivial operator commuting
with the action of $GL(n,F)$. Recall that an irreducible
$GL(n,F)$- module is called {\itshape unramified} if it has a
non-zero vector invariant with respect to maximal compact subgroup
of $GL(n,F)$.

\begin{theorem}\label{T:main-p-adic}
The source and target of operator $T$ are irreducible $GL(n,F)$-representations. Consequently, $T$ is an isomorphism.
\end{theorem}

The proof of this theorem is an application of the results of the
paper by Zelevinsky \cite{zelevinsky}. Let us remind the necessary
facts following the notation of \cite {zelevinsky}.

We will denote for brevity the group $GL(n,F)$ by $G_n$. By $Irr(G_n)$
we will denote the set of isomorphism classes of irreducible
representations of $G_n$, and by $Rep(G_n)$ we will denote the set of isomorphism
classes of all smooth representations of $G_n$.
 Let $\alpha =(n_1, \dots , n_r)$ be an ordered partition
of $n$. Let $G_{\alpha}$ be the subgroup $G_{n_1} \times \dots \times G_{n_r}$
of $G_n$  consisting of block- diagonal matrices. Let $P_{\alpha}$
denote the subgroup of $G_n$ consisting of block-upper-diagonal matrices.
Then $P_{\alpha}$ is a parabolic subgroup with the Levi factor isomorphic to $G_{\alpha}$.
For $\rho _i  \in Rep(G_{n_i}) , \, i= 1, \dots ,r$, let $\rho_1 \otimes \dots \otimes \rho _r \in
Rep(G_{\alpha})$ be their (exterior) tensor product. This representation can be extended
(trivially) to the representation of $P_{\alpha}$. Define
$$ \rho_1 \times \dots \times \rho _r :=
Ind _{P_{\alpha}} ^{G_n}( \rho_1 \otimes \dots \otimes \rho _r), $$
where the induction is normalized.

Let ${\cal C}$ be the set of equivalence classes of irreducible cuspidal
representations of $G_n, \, n=1,2, \dots$. Note that a (one- dimensional) character
of $F^*$ can be considered as a cuspidal representation of $G_1=F^*$.
Let us call a segment in ${\cal C}$ any subset of ${\cal C}$ of the form
$\Delta =\{\rho, \nu \rho , \dots , \nu ^k \rho =\rho '\}$, where $\nu$ is the
character $\nu (g)= |det(g)|$. We will write it also as $\Delta = [\rho , \rho ']$.
To each segment $\Delta$ one can
associate the irreducible representation $<\Delta>$ which can be defined
as the unique irreducible submodule of $\rho \times \nu \rho \times \dots \times \nu ^k \rho$.

Let $\Delta _1 =[\rho _1, \rho _1 ']$ and $\Delta _2=[\rho _2, \rho _2 '] $
be two segments in ${\cal C}$. The segments
$\Delta _1$ and $\Delta _2$ are called {\itshape linked} if $\Delta _1 \cup \Delta _2$
is also a segment and $\Delta _1 \not \subset \Delta _2, \,
\Delta _2 \not \subset \Delta _1$. The segments
$\Delta _1$ and $\Delta _2$ are called {\itshape juxtaposed} if
they are linked and $\Delta _1 \cap \Delta _2 =\emptyset$.
We say that $\Delta _1$ {\itshape precedes} $\Delta _2$
if they are linked and $\rho _2 =\nu ^k \rho _1$ for $k>0$.

The following theorem is one of the main results of \cite{zelevinsky} (Theorem 6.1).
\begin{proposition}
(a) Let $\Delta _1, \dots , \Delta _r$ be segments in ${\cal C}$. Suppose that
for each pair of indices $i<j$, $\Delta _i$ does not precede $\Delta _j$.
Then the representation $<\Delta _1> \times \dots \times <\Delta _r>$
has a unique irreducible submodule; denote it by $<\Delta _1, \dots , \Delta _r>$.

(b) The representations $<\Delta _1, \dots , \Delta _r>$ and $<\Delta '_1, \dots , \Delta '_s>$
are isomorphic iff the sequences
$<\Delta _1>, \dots , <\Delta _r>$ and $<\Delta '_1>, \dots , <\Delta '_s>$
are equal up to rearrangement.

(c) Any irreducible representation of $G_n$ is isomorphic to some representation of the form
$<\Delta _1, \dots , \Delta _r>$.
\end{proposition}

We will also need the following fact (\cite{zelevinsky}, Theorem 4.2).
\begin{proposition}\label{P:not-linked}
Let $\Delta _1, \dots , \Delta _r$ be segments in ${\cal C}$.
Then the representation $<\Delta _1> \times \dots \times <\Delta _r>$ is irreducible
iff for all $i \ne j$ the segments $ \Delta _i$ and $\Delta _j$ are not linked.
\end{proposition}

{\bf Proof  of Theorem \ref{T:main-p-adic}.} We may and will assume that $1\leq i\leq n-1$.  We use the notation of the beginning
of this section. Recall that we study the intertwining operator
$$ T: \Gamma (Gr_{n-i,n}, N) \str \Gamma (Gr_{i,n}, L) .$$
It is easy to see that $\Gamma (Gr_{i,n}, L)$ as a the representation of
$GL(n,F)$ is isomorphic to $<\Delta _1> \times <\Delta _1'>$, where
$$\Delta _1 = (\underset{i\mbox{ times }}{\underbrace{\nu^{-\frac{n+1}{2}},\nu^{-\frac{n-1}{2}},\dots,\nu^{-\frac{n+3}{2}+i}}}),$$
$$ \Delta _1 ' = (\underset{n-i\mbox{ times}}{\underbrace{\nu^{-\frac{n-1}{2}+i},\nu^{-\frac{n-3}{2}+i},\dots,\nu^{\frac{n-1}{2}}}}).$$ 
We see that the segments
$\Delta _1$ and $\Delta _1 '$ are not linked.  Hence by Proposition \ref{P:not-linked} the
representation $\Gamma (Gr_{i,n}, L)=<\Delta _1> \times <\Delta _1'>$ is irreducible. \qed

Let us prove that $\Gamma (\gr_{n-i,n},N)$ is also irreducible representation.
Recall that $N=M \otimes|\omega|$.

It is easy to see that $\Gamma (Gr_{n-i,n},N)$ as a representation of $GL(n,F)$ is isomorphic to 
$$<\Delta_2>\times <\Delta'_2>,$$
where 
\begin{eqnarray*}
\Delta_2=<\underset{n-i\mbox{ times}}{\underbrace{\nu^{-\frac{n-1}{2}+i},\nu^{-\frac{n-3}{2}+i},\dots, \nu^{\frac{n-1}{2}}}}>,\\
\Delta'_2=<\underset{i\mbox{ times}}{\underbrace{\nu^{-\frac{n+1}{2}},\nu^{-\frac{n-1}{2}},\dots,\nu^{-\frac{n+3}{2}+i}}}>.
\end{eqnarray*}
Obviously, the segments $\Delta_2$ and $\Delta_2'$ are not linked. Hence by Proposition \ref{P:not-linked} the
representation $<\Delta _2> \times <\Delta _2'>$ is irreducible.

Since the operator $T$ commutes with $GL(n,F)$, and since its source and the target are irreducible representations, it follows that $T$ is an isomorphism. \qed

\end{document}